\documentclass[11pt]{amsart}
\usepackage{amsmath,amscd,amssymb,amsfonts,latexsym}
\usepackage[all]{xy}
\usepackage{ulem}

\usepackage[latin1]{inputenc}

\usepackage[T1]{fontenc}

\usepackage{eucal} 
\usepackage{cmmib57} 
\usepackage{amsbsy}
\usepackage{amsthm,alltt,dsfont}

\usepackage[dvips]{graphics}

\usepackage{geometry}
\allowdisplaybreaks

\newcommand{\bC}{{\mathbb C}}

\def\bF{\mathbb F}
\newcommand{\bH}{{\mathbb H}}

\newcommand{\bV}{{\mathbb V}}

\newcommand{\cF}{{\mathcal F}}
\newcommand{\cG}{{\mathcal G}}
\newcommand{\cH}{{\mathcal H}}

\newcommand{\cN}{{\mathcal N}}

\newcommand{\cR}{{\mathcal R}}

\newcommand{\cL}{{\mathcal L}}

\newcommand{\cU}{{\mathcal U}}

\newcommand{\Z}{\mathbb{Z}}  
\newcommand{\Q}{\mathbb{Q}}  
\newcommand{\C}{\mathbb{C}}  
\newcommand{\A}{\mathcal{A}}
\newcommand{\CP}{\mathbb{CP}^{n+1}}  
\newcommand{\TA}{H^{\varepsilon,\rho}} 
\newcommand{\TC}{C^{\varepsilon,\rho}} 
\newcommand{\Ft}{\mathbb{F}[t^{\pm 1}]} 

\newcommand{\Hom}{{\rm{Hom}}}
\newcommand{\Ext}{{\rm{Ext}}}

\newcommand{\End}{{\it End}}
\newcommand{\Perv}{{\it Perv}}

\newcommand{\cok}{\it Coker}

\newcommand{\lra}{\longrightarrow}

\newcommand{\rank}{{\it rank}}

\def\eps{\varepsilon}

\theoremstyle{plain}
\newtheorem{thm}{Theorem}[section]
\newtheorem{cor}[thm]{Corollary}

\newtheorem{lem}[thm]{Lemma}
\newtheorem{prop}[thm]{Proposition}

\theoremstyle{definition}
\newtheorem{df}[thm]{Definition}
\newtheorem{rem}[thm]{Remark}
\newtheorem{example}[thm]{Example}

\def\be{\begin{equation}}
\def\ee{\end{equation}}

\def\bt{\begin{thm}}
\def\et{\end{thm}}

\def\bc{\begin{cor}}
\def\ec{\end{cor}}

\def\br{\begin{rem}}
\def\er{\end{rem}}

\def\bp{\begin{prop}}
\def\ep{\end{prop}}

\def\bl{\begin{lem}}
\def\el{\end{lem}}

\def\bn{\begin{enumerate}}
\def\en{\end{enumerate}}

\def\bex{\begin{example}}
\def\eex{\end{example}}

\def\bd{\begin{df}}
\def\ed{\end{df}}

\title[Twisted Alexander invariants]{Twisted Alexander invariants \\ of complex hypersurface complements}
\author{Lauren\c{t}iu Maxim}
\address{L. Maxim: Department of Mathematics, University of Wisconsin-Madison,  480 Lincoln Drive, Madison, WI 53706, 
USA.}
\email {maxim@math.wisc.edu}
\author{Kaiho Tommy Wong}
\address{K. Wong: Department of Mathematics, University of Wisconsin-Madison,  480 Lincoln Drive, Madison, WI 53706, 
USA.}
\email {wong@math.wisc.edu}

\date{\today}



\begin{document}

\begin{abstract}  
We define and study twisted Alexander-type invariants of complex hypersurface complements. 
We investigate torsion properties for the twisted Alexander modules and extend local-to-global divisibility results of \cite{M06,DL} to the twisted setting. In the process, we also study the splitting fields containing the roots of the corresponding twisted Alexander polynomials.
\end{abstract}

\maketitle

\tableofcontents


\section{Introduction}
The classical Alexander polynomial from knot theory has proved to be a powerful and versatile tool in the study of complements of plane algebraic curves. As showed by Libgober in \cite{L82}, the Alexander polynomial of a plane curve complement is sensitive to the local type and position of singularities of the curve, and it can be used to detect Zariski pairs (i.e., pairs of plane curves which have homeomorphic tubular neighborhoods, but non-homeomorphic complements). The study of Alexander polynomials of complements of higher-dimensional complex hypersurfaces have been initiated by Libgober in \cite{L94}, and was pursued in greater generality (for arbitrary singularities) in \cite{M06,DL,Liu}.

A twisted version of the Alexander polynomial (based on the extra datum of a representation of the fundamental group) was introduced by Lin \cite{L}, Wada \cite{W}, Kirk-Livingston \cite{KL} in the $1990$s, and has well proved its worth, for instance, in the works of Friedl and  Vidussi (e.g., see \cite{FV} and the references therein). 

The twisted Alexander polynomial was ported to the study of plane algebraic curves by Cogolludo and Florens \cite{CF}, who used it to refine Libgober's divisibility results from \cite{L82}, and showed that these twisted Alexander polynomials can detect Zariski pairs which were undistinguishable by the classical Alexander polynomial.

In this paper, we extend the Cogolludo-Florens construction to high dimensions and arbitrary singularities, and establish some of the basic properties of the twisted Alexander invariants in this algebro-geometric setting. More concretely, we investigate torsion properties for the twisted Alexander modules, and extend local-to-global divisibility results of \cite{M06,DL} to the twisted setting. In the process, we also study the splitting fields containing the roots of the corresponding twisted Alexander polynomials.


\subsection*{Main results} In what follows, we give a brief overview of our results.

Let $V \subset \CP$ be a projective complex hypersurface, and fix a hyperplane $H$ in $\CP$, which we call the hyperplane at infinity. Let $$\cU:=\CP\setminus (V \cup H)$$ denote the (affine) hypersurface complement. Fix a field $\bF$ which is a subfield of $\C$ closed under conjugation, and let $\bV$ be a finite dimensional $\bF$-vector space. To a pair $(\eps,\rho)$ of an epimorphism $\eps:\pi_1(\cU) \to \Z$ and a representation $\rho:\pi_1(\cU) \to GL(\bV)$, we associate {\it (co)homological (global) twisted Alexander modules} $\TA_i(\cU,\Ft)$ and resp. $H^i_{\eps,\rho}(\cU,\Ft)$, which are $\Ft$-modules of finite type (and which are related by the universal coefficient theorem). These are  homotopy invariants of the complement $\cU$.

\medskip

We say that the hypersurface $V$ is {\it in general position at infinity} if the reduced variety $V_{red}$ underlying $V$ is transversal to $H$ in the stratified sense. One of our first results describes torsion properties of the (global) twisted Alexander modules (see Theorems \ref{t1} and \ref{t3} and Corollarly \ref{imc}):
\bt\label{t1i}
Let $V \subset \CP$ be a hypersurface in general position at infinity.  The twisted Alexander modules $\TA_i(\cU,\Ft)$ are torsion $\Ft$-modules for any $0\leq i \leq n$, they are trivial for $i>n+1$, and $H^{\eps,\rho}_{n+1}(\cU,\Ft)$ is a free $\Ft$-module of rank 
$(-1)^{n+1}\cdot \dim_{\bF}(\bV)\cdot \chi(\cU).$
\et
This is a far-reaching generalization of results from \cite{M06,DL,Liu}, which only dealt with the  case of the linking number homomorphism and the trivial representation defined on complements of {\it reduced} hypersurfaces.

\medskip

For any point $x \in V$, let $\cU_x=\cU \cap B_x$ denote the local complement at $x$, for $B_x$ a small ball about $x$ in $\CP$. Then $(\eps,\rho)$ induces via the inclusion map $i_x:\cU_x \hookrightarrow \cU$ a pair $(\eps_x,\rho_x)$ on $\cU_x$, so that {\it local twisted Alexander modules} of $(\cU_x, \eps_x,\rho_x)$ can be defined. Proposition \ref{ploc} asserts that for any pair $(\eps,\rho)$ as above, we have the following local torsion property:
\bp\label{p1i} If $V$ is in general position at infinity, then the local twisted Alexander modules $H_i^{\eps_x,\rho_x}(\cU_x,\Ft)$ are torsion $\Ft$-modules for any $x \in V$.
\ep 
This local torsion property removes a technical assumption used by Cogolludo-Florens  \cite{CF} in the proof of their main divisibility result for twisted Alexander polynomials of plane curve complements.

\medskip

Since $\Ft$ is a PID, torsion $\Ft$-modules of finite type have orders (called Alexander polynomials) associated to them. Let $\Delta_{i,\cU}$ (resp. $\Delta^i_{\cU}$) and  $\Delta_{i,x}$ (resp. $\Delta^i_x$) be the corresponding {\it global} and {\it local} twisted Alexander polynomial associated to the above (co)homological twisted Alexander modules. In Theorem \ref{t5} we indicate how to compute the global twisted Alexander polynomials from the local topological information at points on the hypersurface. This relationship can be roughly formulated as follows (see Theorem \ref{t5} for the precise folrmulation):
\bt\label{t2i}  For a projective hypersurface $V$ in general position at infinity, the zeros of the global twisted Alexander polynomials of the complement $\cU$ are among those of the local ones at points in the affine part of some irreducible component of $V$. \et
This result is a generalization to the twisted setting of the local-to-global analysis for the classical Alexander polynomials initiated in the first author's work \cite{M06}, and continued in \cite{DL, Liu}.
As a consequence of the proof of Theorem \ref{t5}, we remark that the local torsion property for the twisted Alexander modules at points in $V \cap H$ is enough to conclude that the global twisted Alexander modules are torsion $\Ft$-modules in the desired range. For hypersurfaces in general position at infinity, the local torsion property at points in $V \cap H$  is a consequence of transversality and the K\"unneth formula, see the proof of Proposition \ref{ploc}, but there may be other instances (e.g., for various choices of $(\eps,\rho)$) when it is satisfied.

\medskip

We also single out the contribution of the meridian at infinity (i.e., a meridian loop about $H$) to the global twisted Alexander polynomials, see Theorem \ref{t4} for the precise formulation. For the case of the linking number homomorphism and trivial representation, Theorem \ref{t4} reduces to the fact that the zeros of the classical Alexander polynomials of $\cU$ are roots of unity of order $d=\deg(V)$, a fact shown in \cite{M06,DL} for reduced hypersurfaces.

In the case of reduced plane curves and for $\eps$ the linking number homomorphism, we identify explicitly the splitting fields containing the roots of the corresponding global twisted Alexander polynomials. Similar results were obtained by Libgober \cite{L09} by Hodge-theoretic methods. More precisely, in Theorem \ref{t2} we prove the following:
\bt\label{t13} Let $C$ be a reduced curve of degree $d$ and in general position at infinity. Denote by $x_0$ the (homotopy class of the) meridian about the line $H$ at infinity.
Suppose $\mathbb{F} = \C$, and denote the eigenvalues of $\rho(x_0)^{-1}$ by $\lambda_1,\cdots,\lambda_{\ell}$. Then the roots of $\Delta_{1,\cU}(t)$ lie in the splitting field $\mathbb{S}$ of $\prod_{i=1}^{\ell} (t^d-\lambda_i)$ over $\Q$, which is cyclotomic over $\mathbb{K}=\Q(\lambda_1,\cdots,\lambda_{\ell})$.\et

This result is based on our calculation of the twisted Alexander polynomial for the Hopf link on $d$ components (see Proposition \ref{p1}), which in our geometric situation can be identified with the link of $C$ ``at infinity''.


\subsection*{Acknowledgement.}
This paper was written while the first author visited the Max-Planck-Institut f\"ur Mathematik in  Bonn, and the Institute of Mathematical Sciences at
the Chinese University of Hong Kong. He thanks these institutes for their hospitality and for providing him with excellent working conditions. 
L. Maxim was partially supported by grants from NSF, NSA, by a fellowship from the Max-Planck-Institut f\"ur Mathematik,  Bonn, and by the Romanian Ministry of National Education, CNCS-UEFISCDI, grant PN-II-ID-PCE-2012-4-0156. 
K. Wong gratefully acknowledges the support provided by the NSF-RTG grant \#1502553 at the University of Wisconsin-Madison. 

\section{Twisted chain complexes. Twisted Alexander invariants}\label{secex}

\subsection{Definitions}\label{def}

In this section, we recall the definitions of twisted chain complexes, twisted Alexander modules, and twisted Alexander polynomials of path-connected finite CW-complexes. For more details, see \cite{KL,CF}.

Let $X$ be a path-connected finite CW-complex, with $\pi=\pi_1(X)$, and fix a field $\mathbb{F}$ which is a subfield of $\C$ closed under conjugation. Fix a group homomorphism $$\varepsilon: \pi_1(X) \rightarrow \Z,$$ and note that $\varepsilon$ extends to an algebra homomorphism $$\varepsilon: \mathbb{F}[\pi] \rightarrow \mathbb{F}[\Z] \cong \Ft.$$
Consider a finite dimensional $\mathbb{F}$-vector space $\mathbb{V}$ and a linear representation $$\rho: \pi \rightarrow GL(\mathbb{V}).$$ For simplicity, this representation will also be denoted by  $\bV_{\rho}$.

Let $\widetilde{X}$ be the universal cover of $X$. The cellular chain complex $C_*(\widetilde{X},\mathbb{F})$ of $\widetilde{X}$ is a complex of free left $\mathbb{F}[\pi]$-modules, generated by lifts of the cells of $X$. 
For notational convenience, we follow \cite{KL} and regard $\mathbb{V}$ as a right $\mathbb{F}[\pi]$-module, i.e., with the right $\pi$-action for $v\in \mathbb{V}$ and $\alpha\in \pi$ given by: $$v\cdot \alpha = \rho(\alpha)(v).$$
Also consider the right $\mathbb{F}[\pi]$-module $\Ft \otimes_{\bF} \mathbb{V}$, with $\mathbb{F}[\pi]$-multiplication  induced by $\varepsilon\otimes \rho$ as:
$$(p\otimes v)\cdot \alpha = pt^{\varepsilon(\alpha)} \otimes v\cdot\alpha= pt^{\varepsilon(\alpha)} \otimes \rho(\alpha)v, \ \alpha \in \pi.$$
Let the chain complex of $(X,\varepsilon,\rho)$ be defined as the complex of $\Ft$-modules:
$$\TC_*(X,\Ft) := (\Ft \otimes_{\bF} \mathbb{V}) \otimes_{\mathbb{F}[\pi]}C_*(\widetilde{X},\mathbb{F}),$$ where the $\Ft$-action is given by $$t^n ((p\otimes v)\cdot c) = (t^n \cdot p \otimes v)\cdot c.$$
It is complex of free $\Ft$-modules.

\bd
The {\it $i$-th homological twisted Alexander module} $\TA_i(X,\Ft)$ of the triple $(X, \varepsilon, \rho)$ is the $\Ft$-module defined by:
 $$\TA_i(X,\Ft):=H_i\big(\TC_*(X,\Ft)\big).$$
 Similarly, the {\it $i$-th cohomological twisted Alexander module} $H^i_{\eps,\rho}(X,\Ft)$ of $(X, \varepsilon, \rho)$ is the $\Ft$-module given by:
 $$H^i_{\eps,\rho}(X,\Ft):=H^i\big(\Hom_{\Ft}(C_*^{\eps,\rho}(X,\Ft),\Ft)\big).$$
\ed

\br The classical Alexander modules correspond to the case of the trivial representation $\rho=triv$, i.e., $\bV=\bF=\Q$ and $\rho(x)=1$ for all $x \in \pi$.
\er

The twisted Alexander modules are homotopy invariants. 

The universal coefficient theorem (UCT) applied to the principal ideal domain $\Ft$ yields that:
\be\label{uct}
H^i_{\eps,\rho}(X,\Ft) \cong \Hom_{\Ft}\big(\TA_i(X,\Ft),\Ft\big) \oplus \Ext_{\Ft}\big(\TA_{i-1}(X,\Ft),\Ft\big).
\ee
So, if $\TA_i(X,\Ft)$ are torsion $\Ft$-modules for all $i \leq n$, then $H^i_{\eps,\rho}(X,\Ft)$ are also torsion in the same range.

\medskip

An  equivalent definition of the twisted chain complex of $(X,\varepsilon,\rho)$ was given in \cite{KL}. Let $X_{\infty}$ be the infinite cyclic cover of $X$ associated to $\pi'=\ker\varepsilon$. The chain complex $$C_*(X_{\infty},\mathbb{V}_{\rho}):=\mathbb{V}\otimes_{\mathbb{F}[\pi']} C_*(\widetilde{X}),$$ defined via the restricted actions to $\pi'$, can be regarded as a complex of $\Ft$-modules via the action $t^n \cdot (v\otimes c) = v\cdot\gamma^{-n} \otimes \gamma^n c$, where $\gamma$ is an element in $\pi$ such that $\varepsilon(\gamma)=1$. Then \cite[Theorem 2.1]{KL} states that $C_*(X_{\infty},\mathbb{V}_{\rho})$ and $\TC_*(X,\Ft)$ are isomorphic as a $\Ft$-modules.

\bd\label{ac}
Denote by $\mathbb{F}(t)$ the field of fractions of $\Ft$, and define $$\TC_*(X,\mathbb{F}(t)) = \TC_*(X,\Ft) \otimes \mathbb{F}(t).$$ 
We say that $(X,\varepsilon,\rho)$ is {\it acyclic} if the chain complex $\TC_*(X,\mathbb{F}(t))$ is acyclic over $\mathbb{F}(t)$.\ed
\br Since $\Ft$ is a principal ideal domain, $\mathbb{F}(t)$ is flat over $\Ft$. So, $(X,\varepsilon,\rho)$ is acyclic if and only if $\TA_*(X,\Ft)$ are torsion $\Ft$-modules. \er

Since $\Ft$ is a principal ideal domain and $\bV$ is finite dimensional over $\bF$, the twisted Alexander modules $\TA_*(X,\Ft)$ are finitely generated modules over $\Ft$. Thus they have a direct sum decomposition into cyclic modules. Similar considerations apply for the cohomological invariants.

\bd\label{tap}
The order of the torsion part of $\TA_i(X,\Ft)$ is called the {\it $i$-th homological twisted Alexander polynomial} of $(X,\varepsilon,\rho)$, and is denoted by $\Delta_{i,X}^{\varepsilon,\rho}(t)$. Similarly, we define the {\it $i$-th cohomological twisted Alexander polynomial} of $(X,\varepsilon,\rho)$ to be the order $\Delta^{i}_{\varepsilon,\rho,X}(t)$ of the torsion part of the $\Ft$-module $H^i_{\eps,\rho}(X,\Ft)$.
\ed 

The twisted Alexander polynomials are well-defined up to units in $\Ft$.
Moreover, it follows from (\ref{uct}) that $$\Delta^{i}_{\varepsilon,\rho,X}(t) = \Delta_{i-1,X}^{\varepsilon,\rho}(t).$$

For further use, we also recall here the following fact:

\begin{prop}\label{p0}\cite{KL} If $\varepsilon$ is non-trivial, then $\TA_0(X,\Ft)$ is a torsion $\Ft$-module.
\end{prop} 


\subsection{Examples}\label{ex} 
In this section, we compute the twisted Alexander invariants on several examples with geometric significance. 

\subsubsection{Hopf link with $d$ components.}\label{Hopf}  This example has important consequences in the study of twisted Alexander invariants of plane curve complements. More precisely, for a degree $d$ plane curve $C$ with regular behavior at infinity, the Hopf link with $d$ components is what we call ``the link of $C$ at infinity''.

Recall that a link in $S^3$ is an embedding of a disjoint union of circles (link components) into $S^3$. Throughout this section, let $K$ be the Hopf link with $d$ components in $S^3$, that is, the link with $d\geq 2$ components with the property that the linking number of any two of its components is $1$. 

\bl\label{l1} If $K \subset S^3$ is the Hopf link with $d$ components, then
\be\label{eq1} \pi_1(S^3\setminus K)\cong  \Z \times F_{d-1} \cong \langle x_0,x_1,\cdots,x_{d-1} \vert \  x_0x_ix_0^{-1}x_i^{-1}, i=1,\cdots,d-1 \rangle,\ee
with $F_{d-1}$ the free group on $d-1$ generators.
\el

\begin{proof}
First note that $S^3\setminus K$ is homotopy equivalent to the link exterior associated to the singularity  $\{x^d=y^d\} \subset \C^2$. Equivalently, if  $\A = \{x^d=y^d\}$ is the central line arrangement of $d$ lines in $\C^2$, then $S^3\setminus K \simeq \C^2 \setminus \A$.

On the other hand, it can be easily seen that 
$$\C^2 \setminus \A \simeq  \C^*\times (\mathbb{CP}^1\setminus \{\text{d points}\}).$$
Indeed, the Hopf fibration $\C^2 \setminus \{0\} \to \mathbb{CP}^1$ restricts to a $\C^*$-locally trivial fibration $\C^2 \setminus \A \to \mathbb{CP}^1\setminus \{\text{d points}\}$. Moreover, the latter fibration is trivial, since it can be seen as a restriction of the trivial fibration
$\C^2 \setminus H \to  \mathbb{CP}^1\setminus \{\text{1 point}\}=\C$ obtained from the Hopf fibration by first restricting to the complement of only one line $H$ of $\A$.

Altogether, $$S^3\setminus K \simeq \C^2 \setminus \A \simeq S^1 \times (\bigvee_{d-1} S^1),$$ which yields the desired presentation for $\pi_1(S^3\setminus K)$.
\end{proof}

\br An equivalent presentation of $\pi_1(S^3\setminus K)$ can be obtained by using the van Kampen theorem (e.g., see \cite[Theorem 4.2.17, Proposition 4.2.21]{Di1} and the references therein). More precisely, $\pi_1(S^3\setminus K)$ is called $G(d,d)$ in loc.cit., and has the presentation:
$$\pi_1(S^3\setminus K) \cong  \langle x_0,x_1,\cdots,x_d \ \vert \ x_dx_{d-1}\cdots x_1 x_0^{-1}, x_0x_ix_0^{-1}x_i^{-1}, i=1,\cdots,d \rangle,$$
where the generators $x_1,\cdots,x_d$ correspond to meridian loops about the $d$ lines of $\A$.
\er

We can now compute the twisted Alexander invariants of $S^3 \setminus K$: 
\bp\label{p1} Let $K \subset S^3$ be the Hopf link with $d$ components.
Let $$\varepsilon: \pi_1(S^3\setminus K) \lra \Z$$ be an epimorphism with $$\varepsilon(x_0) \neq 0,$$ and  $$\rho:\pi_1(S^3\setminus K) \lra GL(\mathbb{V}) = GL_{\ell}(\mathbb{F})$$ be a linear representation of rank $\ell$.
Then the following hold:
\begin{itemize}
\item[(a)] $\TA_i(S^3\setminus K,\Ft)$ are torsion $\Ft$-modules, for $i=0,1$.
\item[(b)] $\TA_i(S^3\setminus K,\Ft)=0$ for $i\geq 2$. 
\item[(c)] $\Delta_0^{\varepsilon,\rho}$ is the greatest common divisor of the $\ell \times \ell$ minors of the column matrix  $$\left( \rho(x_i)t^{\varepsilon(x_i)} -Id\right)_{i=0,\cdots,d-1}.$$ 
\item[(d)] $\Delta_1^{\varepsilon,\rho}/\Delta_0^{\varepsilon,\rho} = \big(\det(\rho(x_0)t^{\varepsilon(x_0)}-Id)\big)^{d-2}$.
\end{itemize}
\ep

\begin{proof}
Recall from Lemma \ref{l1} that the link complement $S^3\setminus K$ has the homotopy type of a (central) line arrangement complement, namely $S^3\setminus K \simeq \C^2 \setminus \A$. As such, it has a minimal cell structure (i.e., so that the number of $i$-cells equals its $i$-th Betti number $b_i$, for all $i\geq 0$). Moreover, since $\C^2 \setminus \A$ has the homotopy type of a finite real $2$-dimensional CW-complex, 
it follows that $\TA_i(S^3\setminus K,\Ft)=0$ for $i\geq 3$.

We next note that $S^3\setminus K$ is a $K(\pi,1)$-space, since $\C^2 \setminus \A$ is so, with $\pi=\pi_1(S^3\setminus K)$. Indeed, since $\A$ is defined by a (weighted) homogeneous polynomial, there is a global Milnor fibration $$F \hookrightarrow \C^2 \setminus \A \lra \C^*$$ whose fiber $F$ has the homotopy type of a join of circles. The long exact sequence of homotopy groups for this fibration then yields that $\pi_i( \C^2 \setminus \A)=0$ for all $i \geq 2$. 

Since $S^3\setminus K$ is a $K(\pi,1)$-space, its (twisted) homology can be computed from its (twisted) group homology using Fox calculus (this was the starting point for Wada's construction of twisted Alexander invariants \cite{W}).
So the twisted chain complex of $S^3 \setminus K$ can be identified with the complex of Fox derivatives for the presentation 
$$\pi_1(S^3 \setminus K) \cong \langle x_0,x_1,\cdots,x_{d-1} \vert \  x_0x_ix_0^{-1}x_i^{-1}, i=1,\cdots,d-1 \rangle$$ of Lemma \ref{l1}, and it has the form:
$$0\lra \Ft^{\ell(d-1)} \overset{\partial_2}{\lra} \Ft^{\ell d} \overset{\partial_1}{\lra} \Ft^\ell \lra 0.$$
In particular, as in \cite[Section 4]{KL}, we have that $\partial_1$ is the column matrix with $i$-th entry given by $$\rho(x_i)t^{\varepsilon(x_i)} - Id,$$ which yields the desired description of $\Delta_0^{\varepsilon,\rho}$. Similarly, 
$\partial_2$ is a $(d-1)\times d$ matrix with entries in $M_{\ell}(\Ft)$ given by the matrix of Fox derivatives of the relations, tensored with $\Ft^{\ell}$. Therefore, $\partial_2$ equals
{\footnotesize{\[ \left( \begin{array}{ccccc}
Id-\rho(x_1)t^{\varepsilon(x_1)} & \rho(x_0)t^{\varepsilon(x_0)}-Id & 0 & \cdots & 0 \\
Id-\rho(x_2)t^{\varepsilon(x_2)} & 0 & \rho(x_01)t^{\varepsilon(x_0)}-Id & \cdots & 0 \\
\vdots & \vdots& \vdots & \vdots & \vdots \\
Id-\rho(x_{d-2})t^{\varepsilon(x_{d-2})} & 0 & \cdots & \rho(x_0)t^{\varepsilon(x_0)}-Id & 0 \\
Id-\rho(x_{d-1})t^{\varepsilon(x_{d-1})} & 0 & \cdots & 0 & \rho(x_0)t^{\varepsilon(x_0)}-Id \end{array} \right)\]
}}
Since, by our assumption, 
$\varepsilon(x_0)\neq 0$, this yields that $\ker(\partial_2)=0$. Therefore, $$\TA_2(S^3\setminus K,\Ft)=0.$$ Also, since $\varepsilon$ is non-trivial, we get by Proposition \ref{p0} that $\TA_0(S^3\setminus K,\Ft)$ is a torsion $\Ft$-module.
So, by using the fact that $$\chi(S^3\setminus K)= b_0-b_1+b_2 = 1-d+(d-1)= 0,$$ we obtain that $$\rank_{\Ft}(\TA_1(S^3\setminus K,\Ft)) = -\chi(S^3\setminus K) =0.$$ Hence the first twisted Alexander module $\TA_1(S^3\setminus K,\Ft)$ is also torsion over $\Ft$. Finally, by \cite[Theorem 4.1]{KL}, we get that $$\Delta_1^{\varepsilon,\rho}/\Delta_0^{\varepsilon,\rho} = \big(\det(\rho(x_0)t^{\varepsilon(x_0)}-Id)\big)^{d-2}.$$
\end{proof}

\subsubsection{Links of $A_{odd}$-singularities}
Let $C=\{ x^2-y^{2n}=0\} \subset \C^2$, and fix $(\varepsilon,\rho)$ as before, with $\varepsilon$ non-trivial. The germ $(C,0)$ of $C$ at the origin of $\C^2$ is known as the ${A}_{2n-1}$-singularity. The curve $C$ is the union of two smooth curves which intersect non-transversely at the origin. Let $K \subset S^3$ be the link of $(C,0)$. Since the defining polynomial of $(C,0)$ is weighted homogeneous, it follows that $S^3 \setminus K \simeq \C^2 \setminus C$ fibers over $S^1 \simeq \C^*$, with fiber homotopy equivalent to a join of circles. In particular, $S^3 \setminus K$ is aspherical, so its twisted Alexander invariants can be computed by Fox calculus from a presentation of the fundamental group.
By \cite{O1}, we have that 
$$\pi_1(S^3 \setminus K)\cong \pi_1(\C^2\setminus C) \cong G(2,2n) = \langle a_i,\beta \ \vert \ \beta = a_1a_0, R_1, R_2 \rangle,$$ where $$R_1: a_{i+2n} =a_i, i=0,...,2n-1, \ \ {\rm and} \ \ \ R_2: a_{i+2} = \beta^{-1}a_i\beta , i=0,...,2n-1.$$
So, explicitly, $$\pi_1(S^3 \setminus K)\cong
\begin{array}{c}
\langle a_0,a_1,...,a_{2n-1},\beta \ \vert \ a_1a_0\beta^{-1},\\
\beta a_2 \beta^{-1}a_0^{-1}, \beta a_4 \beta^{-1}a_2^{-1},...,\beta a_0 \beta^{-1} a_{2n-2}^{-1},\\
\beta a_3 \beta^{-1}a_1^{-1}, \beta a_5 \beta^{-1}a_3^{-1},...,\beta a_1 \beta^{-1} a_{2n-1}^{-1}\rangle
\end{array}$$
By direct computation, it can be seen that in the corresponding twisted chain complex one has  $\ker(\partial_2)=0$, so $\TA_2(S^3\setminus K,\Ft)=0.$ Also, since $\varepsilon$ is non-trivial, we get by Proposition \ref{p0} that $\TA_0(S^3\setminus K,\Ft)$ is a torsion $\Ft$-module. An Euler characteristic argument similar to that of the previous example then yields that $\TA_1(S^3\setminus K,\Ft)$ is a torsion $\Ft$-module.


\section{Twisted Alexander invariants of plane curve complements}\label{tc}
Twisted Alexander invariants were ported to the study of plane algebraic curves by Cogolludo and Florens \cite{CF}, who showed that these twisted invariants can detect Zariski pairs which share the same (classical) Alexander polynomial. In this section, we study torsion properties of the twisted Alexander modules of plane curve complements and study the splitting fields containing the roots of the corresponding twisted Alexander polynomials. We focus here on homological invariants, while similar statements about their cohomological counterparts can be obtained via the universal coefficient theorem (\ref{uct}).

Let $C$ be a reduced curve in $\C \mathbb{P}^2$ of degree $d$ with $r$ irreducible components, and let $L$ be a line in $\C \mathbb{P}^2$. Set $$\cU:=\C \mathbb{P}^2 \setminus (C\cup L) = \C^2 \setminus (C\setminus (C\cap L)),$$ where we use the natural identification of $\C^2$ with $\C \mathbb{P}^2 \setminus L$. The line $L$ will usually be refered to as the {\it line at infinity}. Alternatively, let $f(x,y): \C^2 \rightarrow \C$ be a square-free polynomial of degree $d$ defining an affine plane curve $C^a:=\{f=0\}$. Let $C$  be the zero locus in $\C \mathbb{P}^2$ (with homogeneous coordinates $x,y,z$) of the projectivization $\bar f$ of $f$, and let $L$ by given by $z=0$. Then $\cU=\C^2 \setminus C^a$.

Recall that $H_1(\cU,\Z)\cong \Z^r$, generated by homology classes $\nu_i$ of meridian loops $\gamma_i$ bounding transversal disks at a smooth point in each irreducible component of $C^a$. Let $n_1,\cdots,n_r$ be positive integers with gcd($n_1,\cdots,n_r$)=1. Let $ab: \pi_1 (\cU) \to H_1(\cU,\Z)$ denote the abelianization map, sending $[\gamma_i]$ to $\nu_i$. Then the composition $$\varepsilon: \pi_1 (\cU) \xrightarrow{ab} H_1(\cU,\Z) \xrightarrow{\psi: \nu_i \rightarrow n_i} \Z$$ defines an epimorphism. If all $n_i =1$, then $\varepsilon$ can be identified with the total linking number homomorphism $$lk: \pi_1 (\cU) \xrightarrow{[\alpha] \mapsto lk(\alpha,C \cup -dL)} \Z,$$ 
which is just the homomorphism $f_\#:\pi_1(\cU) \to \pi_1(\C^*)\cong \Z$ induced by the restriction of $f$ to $\cU$ (e.g., see \cite[pp.77]{Di1}).

Fix a finite dimensional $\mathbb{F}$-vector space $\mathbb{V}$ and a linear representation $\rho: \pi_1(\cU) \rightarrow GL(\mathbb{V}).$ As in Section \ref{def}, the $\Ft$-module $\TA_i(\cU,\Ft)$ is defined for any $i\geq 0$, and is called the {\it $i$-th (homological) twisted Alexander modules of $C$ with respect to $L$}. 
The twisted Alexander modules associated to the total linking number homomorphism $lk$ will be denoted by $$H_i^{\rho}(\cU,\Ft) := H^{lk,\rho}_i(\cU,\Ft).$$
In the case of the trivial representation, these further reduce to the classical Alexander modules, as originally studied in \cite{L82}.

Note that since $\cU$ is the complement of a plane affine curve, it is a complex $2$-dimensional Stein manifold. Therefore $\cU$ has the homotopy type of a real $2$-dimensional finite CW-complex. Hence, $\TA_i(\cU,\Ft) =0$ for $i\geq3$, and $\TA_2(\cU,\Ft)$ is a free $\Ft$-module. For $i=0,1$, the $\Ft$-modules $\TA_i(\cU,\Ft)$ are of finite type, and in the next section we investigate their torsion properties.


\subsection{Torsion properties}\label{torc}
In this section, we use the above notations to prove the following result:

\bt\label{t1}
Let $C$ be a reduced complex projective plane curve. If $C$ is irreducible and $\rho$ is abelian (i.e., the image of $\rho$ is abelian), or if $C$ is in general position at infinity (i.e., $C$ is transversal to the line at infinity $L$), then the twisted Alexander modules $\TA_i(\cU,\Ft)$ are torsion $\Ft$-modules, for $i=0,1$.
\et

\begin{proof}
The claim about $\TA_0(\cU,\Ft)$ follows from Proposition \ref{p0} since $\varepsilon$ is non-trivial.

If $C$ is irreducible and $\rho$ is abelian, it follows from \cite{L09} that the classical Alexander modules of an irreducible curve complement determine the twisted ones. So the claim follows in this case from \cite{L82}.

Assume now that the line at infinity $L$ is transversal to the curve $C$, and let $d=\deg(C)$. Let $S^3_{\infty} \subset \C^2$ be a sphere of sufficiently large radius. Then the link of $C$ at infinity, $K_{\infty} = S^3_{\infty} \cap C$, is the Hopf link  on $d$ components, as described in Section \ref{Hopf}. Let $i:S^3_{\infty} \setminus K_{\infty} \hookrightarrow \cU$ denote the inclusion map. 
Then by \cite[Lemma 5.2]{L82}, the induced homomorphism  
$$\pi_1(S^3_{\infty} \setminus K_{\infty})\cong \langle x_0,x_1,\cdots ,x_d \ \vert \  x_dx_{d-1}\cdots x_1 x_0^{-1}, x_0x_ix_0^{-1}x_i^{-1}, i=1,\cdots,d\rangle  \overset{i_\#}{\lra} \pi_1(\cU)$$ is surjective. Moreover, as in \cite[Section 7]{L82}, the groups $\pi_1(\cU)$ and $\pi_1(S^3_{\infty} \setminus K_{\infty})$ have the same generators, while the relations in $\pi_1(\cU)$ are those of $\pi_1(S^3_{\infty} \setminus K_{\infty})$ together with relations describing the monodromy about exceptional lines by using the Zariski-Van Kampen method. Therefore, $\varepsilon \circ i_\#= \varepsilon$ and $\rho\circ i_\# = \rho$ (as this can be checked on generators).

Up to homotopy, $\cU$ is obtained from $S^3_{\infty} \setminus K_{\infty}$ by attaching cells of dimension $\geq 2$. So the homomorphism 
$$\TA_k(S^3_{\infty} \setminus K_{\infty},\Ft) \lra \TA_k(\cU,\Ft)$$
induced by the inclusion map $i$ is an isomorphism for $k=0$, and an epimorphism for $k=1$.
Here, $\TA_k(S^3_{\infty} \setminus K_{\infty},\Ft)$ is defined with respect to the pair $(\varepsilon \circ i_\#= \varepsilon, \rho\circ i_\# = \rho)$ induced by the inclusion map $i$.
As a consequence, in order to conclude that $\TA_1(\cU,\Ft)$ is a $\Ft$-torsion module, 
it suffices by Proposition \ref{p1} to show that $\varepsilon\circ  i_\# (x_0) =\varepsilon(x_0)\neq 0$. 

We have a commutative diagram:
$$\xymatrix{
\pi_1(S^3_{\infty} \setminus K_{\infty})\ar[d]^{ab} \ar[r]^{i_\#} & \pi_1(\cU)\ar[r]^{\varepsilon} \ar[d]^{ab} & \mathbb{Z}
\\
H_1(S^3_{\infty} \setminus K_{\infty},\Z) \ar[r]^{i_*} & H_1(\cU,\Z)   \ar[ru]^{\psi}  }$$
So, $\varepsilon\circ i_\# = \psi\circ i_*\circ ab$, hence it is enough to understand the maps $ab$ and $i_*$.
 Recall that the Hopf link complement $S^3_{\infty} \setminus K_{\infty}$ is homotopy equivalent to the complement $\C^2 \setminus \A$ of a central line arrangement $\A$ of $d$ lines in $\C^2$. So 
$$H_1(S^3_{\infty} \setminus K_{\infty},\Z) \cong \Z^d = \langle \mu_1,...,\mu_d \rangle,$$ where $\mu_k$ is the homology class of the meridian about the line $l_k\subset \A$. Moreover,  $ab(x_k) = \mu_k$ for $k=1,\cdots,d$, hence $$ab(x_0) = \mu_1+\cdots+\mu_d.$$
On the other hand, $H_1(\cU,\Z) = \Z^r$, generated by the homology classes $\nu_l$ of the meridians about each irreducible component of $C^a$. Since $\A$ is defined by the homogeneous part of the defining equation of $C^a$, it is clear that $i_*$ takes each $\mu_k$ to one of the $\nu_l$'s. In fact, exactly $d_l$ of the $\mu_k$'s  are being mapped by $i_*$ to $\nu_l$, where $d_l$ is the degree of the component $C_l$ of $C$. 
Finally, since $\psi(\nu_l) =n_l$, for all $k\geq 1$ we have that $\varepsilon\circ i_\#(x_k) =n_{l_j}$ for some $l_j$, and $$\varepsilon\circ i_\#(x_0) = \psi\circ i_* (\mu_1+\cdots+\mu_d) =\sum_{l=1}^r d_ln_l> 0.$$
This concludes the proof of the fact that $\TA_1(\cU,\Ft)$ is a finitely generated $\Ft$-torsion module.
\end{proof}

\br The above result will be generalized in Theorem \ref{t3} to arbitrary (possibly non-reduced) hypersurfaces. The reason for stating it in this section is our study of splitting fields containing the roots of the associated twisted Alexander polynomials, see Theorem \ref{t2}.
\er


As a consequence of Theorem \ref{t1} and  Proposition \ref{p1}, we obtain the following:

\bc\label{c1}
If $C$ is a reduced curve of degree $d$ in general position at infinity, then the first twisted Alexander polynomial $\Delta_{1,\cU}^{\varepsilon,\rho}(t)$ of $\cU$ divides 
{\footnotesize{$$\gcd(\det(\rho(x_0)t^{\sum_{l=1}^r d_ln_l}-Id), \det(\rho(x_1)t^{n_{l_1}}-Id),\cdots,\det(\rho(x_{d-1})t^{n_{l_{d-1}}}-Id))\cdot(\det(\rho(x_0)t^{\sum_{l=1}^r d_ln_l}-Id))^{d-2}.$$}} 
In particular, if $\varepsilon = lk$, then $\Delta_{1,\cU}^{\rho}(t)$ divides $$\gcd(\det(\rho(x_0)t^d-Id),\det(\rho(x_{1})t-Id),\cdots,\det(\rho(x_{d-1})t-Id))\cdot(\det(\rho(x_0)t^d-Id))^{d-2}.$$
\ec

\br For curves in general position at infinity, Corollary \ref{c1} generalizes Libgober's divisbility result \cite[Theorem 2]{L82}, which states that the Alexander polynomial $\Delta_{1,\cU}(t):=\Delta_{1,\cU}^{lk,triv}(t)$ of $C$ divides the Alexander polynomial of the link at infinity, which is given by $(t-1)(t^d-1)^{d-2}$.
\er


\subsection{Roots of twisted Alexander polynomials}\label{rc}

In \cite[Theorem 5.4]{L09}, Libgober used Hodge theory to show that for an irreducible curve $C$, and for $\rho$ a unitary representation, the roots of the first twisted Alexander polynomial of $C$ are in a cyclotomic extension of the field generated by the rationals and the eigenvalues of $\rho(\gamma)$, where $\gamma$ is a meridian about $C$ at a non-singular point. Libgober's result does not touch upon the extension degree. 

In this section, we give a topological proof of Libgober's result, and identify this cyclotomic extension in an explicit way.

\bt\label{t2} Let $C$ be a reduced curve of degree $d$ and assume that $C$ is in general position at infinity. Denote by $x_0$ the (homotopy class of the) meridian about the line $L$ at infinity.
Suppose $\mathbb{F} = \C$, and denote the eigenvalues of $\rho(x_0)^{-1}$ by $\lambda_1,\cdots,\lambda_{\ell}$. Then the roots of $\Delta_{1,\cU}^{\rho}(t)$ lie in the splitting field $\mathbb{S}$ of $\prod_{i=1}^{\ell} (t^d-\lambda_i)$ over $\Q$, which is cyclotomic over $\mathbb{K}=\Q(\lambda_1,\cdots,\lambda_{\ell})$.
\et

\begin{proof} Let us denote as before by $x_1,\cdots,x_d$ the (homotopy classes of) meridians about the irreducible components of $C$. 

If there is no common eigenvalue for all of $\rho(x_1),\cdots,\rho(x_d)$, then Corollary \ref{c1} yields that $\Delta_{1,\cU}^{\rho}(t)$ divides $(\det(\rho(x_0)t^d-Id))^{d-2}.$ In particular, the prime factors of $\Delta_{1,\cU}^{\rho}(t)$ are among the prime factors of $\det(\rho(x_0)t^d-Id)$.
Let $p(t)$ be the characteristic polynomial of $\rho(x_0)^{-1}$.
Then: $$\det(\rho(x_0)t^d-Id) = (-1)^r\det(\rho(x_0))\cdot p(t^d)=(-1)^r\det(\rho(x_0))\cdot (t^d-\lambda_1)\cdots (t^d-\lambda_{\ell}).$$
Therefore,  the roots of $\Delta_{1,\cU}^{\rho}(t)$ are contained in the splitting field $\mathbb{S}$ of $\prod_{i=1}^{\ell} (t^d-\lambda_i)$ over $\Q$.

If $\alpha$ is a common eigenvalue of all matrices $\rho(x_1),\cdots,\rho(x_d)$, then one of the eigenvalues of $\rho(x_0)=\rho(x_d)\rho(x_{d-1})...\rho(x_1)$ is $\alpha^d$. Without loss of generality, assume that $\alpha^d=\lambda_1^{-1}$. Then $\alpha\in \mathbb{S}$.
\end{proof}


\section{Twisted Alexander invariants of complex hypersurface complements}\label{th}
In this section, we generalize the above results to the context of complex hypersurfaces with arbitrary singularities. We study the torsion properties of the associated twisted Alexander modules, and compute their corresponding twisted Alexander polynomials in terms of local topological data encoded by the singularities.


\subsection{Definitions}
Let $V$ be a (not necessarily reduced) degree $d$ hypersurface in $\CP$ ($n \geq 1$) and let $H$ be a hyperplane in $\CP$, called the ``hyperplane at infinity''. Let 
$$\cU:=\CP \setminus (V\cup H) = \C^{n+1} \setminus V^a,$$
where $V^a \subset \C^{n+1} =\CP \setminus H$ denotes the affine part of $V$.
Alternatively, if $f(z_1,\cdots,z_{n+1}): \C^{n+1} \rightarrow \C$ is a polynomial of degree $d$, then $V^a=\{f=0\}$ and $V \subset \CP$ is the projectivization of $V_a$, with $H$ given by $z_0=0$.

Assume that the underlying reduced hypersurface $V_{red}$ of $V$ has $r$ irreducible components $V_1, \cdots, V_r$, with $d_i=\deg(V_i)$ for $i=1,\cdots,r$. Then $$H_1(\cU,\Z) \cong \Z^r,$$  generated by the homology classes $\nu_i$ of meridians $\gamma_i$ about the irreducible components $V_i$ of $V_{red}$ (e.g., see \cite{Di1}, (4.1.3), (4.1.4)). 
Moreover, if $\gamma_{\infty}$ denotes the meridian loop in $\cU$ about the hyperplane $H$ at infinity, with homology class $\nu_{\infty}$, then the following relation holds in $H_1(\cU,\Z)$:
\be\label{relinf} \nu_{\infty}+\sum_{i=1}^r d_i \nu_i=0.\ee  

Let $n_i$ be $r$ positive integers with $\gcd(n_1,\cdots,n_r)=1$, and define the epimorphism $\eps:\pi_1 (\cU) \to \Z$ by the composition $$\varepsilon: \pi_1 (\cU) \xrightarrow{ab} H_1(\cU,\Z) \xrightarrow{\nu_i \mapsto n_i} \Z.$$
Note that if the defining equation $f$ of the affine hypersurface $V^a$ has an irreducible decomposition given by $f=f_1^{n_1}\cdots f_r^{n_r}$, then $\eps$ coincides with 
the homomorphism $f_\#:\pi_1(\cU) \to \pi_1(\C^*)\cong \Z$ induced by the restriction of $f$ to $\cU$, or equivalently, 
with the {\it total linking number homomorphism} (cf. \cite[p.76-77]{Di1}):
$$lk: \pi_1 (\cU) \xrightarrow{[\alpha] \rightarrow lk(\alpha,V \cup -dH)} \Z.$$

Fix a finite dimensional $\mathbb{F}$-vector space $\mathbb{V}$ and a linear representation $\rho: \pi_1(\cU) \rightarrow GL(\mathbb{V}).$ As in Section \ref{def}, the $\Ft$-modules $\TA_i(\cU,\Ft)$ and $H^i_{\eps,\rho}(\cU,\Ft)$ are defined for any $i\geq 0$, and are called the {\it $i$-th (co)homological twisted Alexander modules of $V$ with respect to the hyperplane at infinity $H$}. 
The twisted Alexander modules associated to the total linking number homomorphism $lk$ will be denoted by $$H_i^{\rho}(\cU,\Ft) := H^{lk,\rho}_i(\cU,\Ft),$$ and similarly for their cohomology counterparts $H^i_{\rho}(\cU,\Ft)$.
In the case of the trivial representation, these further reduce to the classical Alexander modules, as studied e.g., in \cite{M06}, \cite{DL} and \cite{Liu}.

Note that, since $\cU$ is the complement of a complex $n$-dimensional affine hypersurface, it  is an $(n+1)$-dimensional affine variety, hence it has the homotopy type of a finite CW-complex of real dimension $n+1$ (e.g., see  \cite{Di1}, (1.6.7), (1.6.8)). 
Therefore, $\TA_i(\cU,\Ft) =0$ for $i\geq n+1$,  $\TA_{n+1}(\cU,\Ft)$ is a free $\Ft$-module,  and the $\Ft$-modules $\TA_i(\cU,\Ft)$ are of finite type for $0\leq i \leq n$. In the next sections, we investigate torsion properties of the latter.


\subsection{Torsion properties}\label{torh}

In the notations of the previous section, we say that the hypersurface $V \subset \CP$ is {\it in general position (with respect to the hyperplane $H$)  at infinity} if the reduced underlying variety $V_{red}$ is transversal to $H$ in the stratified sense.

The main result of this section is the following high-dimensional generalization of Theorem \ref{t1}:
\bt\label{t3}
If the hypersurface $V \subset \CP$ is in general position at infinity, then for any $0\leq i \leq n$ the twisted Alexander modules $\TA_i(\cU,\Ft)$ are torsion $\Ft$-modules.
\et

In order to prove Theorem \ref{t3}, we need to introduce some notations and develop some prerequisites.

\medskip

Let $S^{2n+1}_{\infty}$ be a $(2n+1)$-sphere in $\C^{n+1}$ of a sufficiently large radius (that is, the boundary of a small tubular neighborhood in $\CP$ of the hyperplane $H$ at infinity). Denote by $$K_{\infty}=S^{2n+1}_{\infty} \cap V^a$$ the {\it link of $V^a$ at infinity}, and by  
$$\cU^{\infty}=S^{2n+1}_{\infty} \setminus K_{\infty}$$ its complement in $S^{2n+1}_{\infty}$. Note that $\cU^{\infty}$ is homotopy equivalent to $T(H) \setminus (V \cup H)$, where $T(H)$ is the tubular neighborhood of $H$ in $\CP$ for which $S^{2n+1}_{\infty}$ is the boundary.  Then a classical argument based on the Lefschetz hyperplane theorem yields that the homomorphism $$\pi_i(\cU^{\infty}) \lra \pi_i(\cU)$$ induced by inclusion is an isomorphism for $i < n$ and it is surjective for $i=n$; see \cite[Section 4.1]{DL} for more details.  It follows that 
\be\label{eq}
\pi_i(\cU,\cU^{\infty})=0 \  \ \ {\rm for \ all} \  \ i \leq n,\ee hence $\cU$ has the homotopy type of a CW complex obtained from $\cU^{\infty}$ by adding cells of dimension $\geq n+1$. 

We denote by $(\eps_{\infty},\rho_{\infty})$ the epimorphism and resp. representation on $\pi_1(\cU^{\infty})$ induced by composing $(\eps,\rho)$ with the homomorphism $\pi_1(\cU^{\infty}) \to \pi_1(\cU)$. Hence the {\it twisted Alexander modules of $V$ at infinity}, $H^{\eps_{\infty},\rho_{\infty}}_i(\cU^{\infty},\Ft)$, can be defined (and similarly for the corresponding cohomology modules). Then the fact that twisted Alexander modules are homotopy invariants yields the following:
\bp\label{p3}
The inclusion map $\cU^{\infty} \hookrightarrow \cU$ induces $\Ft$-module isomorphisms 
$$H^{\eps_{\infty},\rho_{\infty}}_i(\cU^{\infty},\Ft) \overset{\cong}{\lra} \TA_i(\cU,\Ft)$$ for any $i<n$, and an epimorphism of $\Ft$-modules
$$H^{\eps_{\infty},\rho_{\infty}}_n(\cU^{\infty},\Ft) \twoheadrightarrow \TA_n(\cU,\Ft).$$
\ep

\bc\label{c2} 
For any $0 \leq i\leq n$, if $H^{\eps_{\infty},\rho_{\infty}}_i(\cU^{\infty},\Ft)$ is a torsion $\Ft$-module, then so is $ \TA_i(\cU,\Ft)$.
\ec

\medskip

Let us now assume that the projective hypersurface $V$ is in general position at infinity, i.e., $V_{red}$ is transversal in the stratified sense to the hyperplane at infinity $H$. 
Then the complement of the link at infinity $\cU^{\infty}$ is a circle fibration over $H \setminus  (V \cap H)$,  which is homotopy equivalent to the complement in $\C^{n+1}$ to the affine cone over the projective hypersurface $V\cap H \subset H=\mathbb{CP}^{n}$ (for a similar argument see \cite[Section~4.1]{DL}). Hence, by the Milnor fibration theorem (e.g., see \cite[(3.1.9),(3.1.11)]{Di1}), $\cU^{\infty}$ fibers over $\C^* \simeq S^1$, with fiber homotopy equivalent to a finite $n$-dimensional CW-complex. Moreover, it is known that this fiber is also homotopy equivalent to the infinite cyclic cover of $\cU^{\infty}$ defined by the kernel of the total linking number homomorphism defined with respect to $V^a$.

We can now complete the proof of Theorem \ref{t3}. 
\begin{proof} By the above Corollary \ref{c2}, it suffices to prove that for any $0\leq i \leq n$, the $\Ft$-module $H^{\eps_{\infty},\rho_{\infty}}_i(\cU^{\infty},\Ft)$ is torsion. The idea is to replace $V^a$ by another affine hypersurface $X$ with the same underlying reduced structure, hence also the same complement $\cU$, so that $\eps$ becomes the  homomorphism defined by the total linking number with $X$.

Let $f_1\cdots f_r=0$ be a square-free polynomial equation defining $V^a_{red}$, the  reduced affine hypersurface underlying $V^a=V \setminus H$. Recall that if $\gamma_i$ is the meridian about the irreducible component $f_i=0$, then by definition we have that $\eps([\gamma_i])=n_i$. Let us now consider the polynomial $g={f_1}^{n_1}\cdots {f_r}^{n_r}$ on $\C^{n+1}$ defining an affine hypersurface $$X=\{g=0\},$$ and replace $V$ by the projective hypersurface $\overline{X}$ defined by the projectivization of $g$. Clearly, the underlying reduced hypersurface $X_{red}$ coincides with $V^a_{red}$, so $X$ and $V^a$ have the same complement $$\cU:=\C^n\setminus V^a=\C^n \setminus X.$$ 
Moreover, the given homomorphism $\eps:\pi_1(\cU)\to \Z$ (hence also $\eps_{\infty}:\pi_1(\cU^{\infty}) \to \Z$) coincides with the total linking number homomorphism defined with respect to $X$  (cf. \cite[p.76-77]{Di1}). Finally, since $V$ is in general position at infinity, so is $\overline{X}$, and the corresponding complements of the links at infinity coincide.
Therefore, the complement of the link at infinity $\cU^{\infty}$ admits a locally trivial topological fibration $$F \hookrightarrow \cU^{\infty} \lra \C^*$$ whose fiber $F$ has the homotopy type of a finite $n$-dimensional CW-complex, and which is also homotopy equivalent to the infinite cyclic cover of $\cU^{\infty}$ defined by the kernel of the linking number with respect to $X$ (i.e., by $\ker(\eps_{\infty})$). 

Altogether, for any $0 \leq i\leq n$, we have:
$$H^{\eps_{\infty},\rho_{\infty}}_i(\cU^{\infty},\Ft) \cong H_i(F, \bV_{\rho_{\infty}}),$$
which is a finite dimensional $\mathbb{F}$-vector space, hence a torsion $\Ft$-module.
\end{proof}

As an immediate consequence of Theorem \ref{t3}, we have the following:
\bc\label{imc} If the hypersurface $V \subset \CP$ is in general position at infinity, then 
$$\rank_{\Ft} H^{\eps,\rho}_{n+1}(\cU,\Ft)=(-1)^{n+1}\cdot \ell \cdot \chi(\cU),$$
with $\ell$ the rank of the representation $\rho$.
\ec

By (\ref{uct}) and Theorem \ref{t3}, we also deduce the following:
\bc\label{cst} If $V \subset \CP$ is a hypersurface in general position at infinity, then the cohomological twisted Alexander modules $H^i_{\eps,\rho}(\cU,\Ft)$ are torsion $\Ft$-modules for any $0\leq i \leq n$.
\ec

\br If $V$ is in general position at infinity, and $\dim_{\C} Sing(V) \leq n-2$ (in which case $V$ is already irreducible), then $\pi_1(\cU) \cong \Z$ (e.g., see \cite[Lemma 1.5]{L94}). So in this case, the representation $\rho$ is abelian, and the twisted Alexander invariants of $\cU$ are determined by the classical ones (already studied in \cite{M06,DL,Liu}). Results of this paper are particularly interesting for hypersurfaces with singularities in codimension one (e.g., hyperplane arrangements) and non-abelian representations.
\er


\subsection{Local twisted Alexander invariants}\label{loca}
For each point $x \in V$, consider the local complement 
$$\cU_x:=\cU \cap B_x,$$ for $B_x$ a small open ball about $x$ in $\CP$ chosen so that $(V,x)$ has a conic structure in $\overline{B}_x$. Let $$\eps_x:\pi_1(\cU_x) \overset{(i_x)_{\#}}{\lra} \pi_1(\cU) \overset{\eps}{\lra}\Z$$ and 
$$\rho_x:\pi_1(\cU_x) \overset{(i_x)_{\#}}{\lra} \pi_1(\cU) \overset{\rho}{\lra}GL(\bV)=GL_{\ell}(\bF)$$
be induced by the inclusion $i_x:\cU_x \hookrightarrow \cU$. The corresponding {\it local} (co)homological twisted Alexander modules $H_k^{\eps_x,\rho_x}(\cU_x,\Ft)$ and $H^k_{\eps_x,\rho_x}(\cU_x,\Ft)$ inherit $\Ft$-module structures.

\br Note that $\eps_x$ is not necessarily onto, so the infinite cyclic cover of $\cU_x$ defined by $\ker(\eps_x)$ may be disconnected. 
\er

\bd\label{ac2} We say that $(\eps,\rho)$ is {\it acyclic at $x \in V$} if $(\eps_x,\rho_x)$ is acyclic in the sense of Definition \ref{ac}, i.e., if $H_k^{\eps_x,\rho_x}(\cU_x,\Ft)$ are torsion $\Ft$-modules for all $k \in \Z$. (Note that by UCT (\ref{uct}), this condition implies that the local cohomological twisted Alexander modules are torsion as well.)
We say that $(\eps,\rho)$ is {\it locally acyclic} along a subset $Y \subseteq V$ if $(\eps,\rho)$ is acyclic at any point $x \in Y$.
\ed

The next result provides one important geometric example of local acyclicity.
\bp\label{ploc} Let $V \subset \CP$ be a degree $d$ projective hypersurface in general position at infinity. Then $(\eps,\rho)$ is locally acyclic along $V$, for any nontrivial epimorphism $\eps:\pi_1(\cU) \to \Z$ and any representation $\rho:\pi_1(\cU)\to GL(\bV)$.\ep

\begin{proof} As in the proof of Theorem \ref{t3}, after changing $V^a$ (resp. $V$) by an affine hypersurface $X$ (resp., by its projectivization $\overline{X}$) with the same underlying reduced structure, hence also preserving the (local) complements, we can assume without loss of generality (and without changing the notations) that $\eps$ is the total linking number homomorphism $lk$. Therefore, for any $x \in V$, the local homomorphism $\eps_x$ becomes $lk_x:=lk \circ (i_x)_{\#}.$ Denote by $\cU_{x,\infty}$ the infinite cyclic cover of $\cU_x$ defined by $\ker(lk_x)$.

Let $\cU'=\CP \setminus V$, and for any point $x\in V$ let $\cU'_x:=\cU' \cap B_x$, for $B_x$ denoting as before a small open ball about $x$ in $\CP$ for which $(V,x)$ has a conic structure in $\overline{B}_x$. Let $S_x:=\partial \overline{B}_x$, with $K_x:=V \cap S_x$ denoting the corresponding link of $(V,x)$. 
Note that $\cU'_x$ is homotopy equivalent to the link complement $S_x \setminus K_x$. Moreover, since $K_x$ is an algebraic link, the Milnor fibration theorem implies that the complement $S_x \setminus K_x$ fibers over a circle, with (Milnor) fiber $F_x$ homotopy equivalent to a finite CW-complex. It is also known that $F_x$ is homotopy equivalent to the infinite cyclic cover of $S_x \setminus K_x$ defined by the linking number with respect to $K_x$. For future reference, let us denote by $lk'_x$ the epimorphism on $\pi_1(S_x \setminus K_x)\cong \pi_1(\cU'_x)$ defined by the total linking number with $K_x$.

If $x \in V \setminus H$, then $\cU_x=\cU'_x \simeq S_x \setminus K_x$, so in this case $$H_k^{\eps_x,\rho_x}(\cU_x,\Ft)=H_k^{lk_x,\rho_x}(\cU_x,\Ft) \cong H_k(\cU_{x,\infty},\bV_{\rho_x}) \cong H_k(F_x,\bV_{\rho_x})$$ is a finite dimensional $\bF$-vector space, hence a torsion $\Ft$-module for any $k \in \Z$. 

If $x \in V \cap H$, then by the transversality assumption we have that $\cU_x \simeq \cU'_x \times S^1$, with the restrictions of $lk_x$ to the factors of this product described as follows: on $\pi_1(\cU'_x)$, $lk_x$ restricts to the homomorphism $lk'_x$ defined by the linking number with $K_x$ (this is, of course, the same as $lk_{x'}$ at a nearby point $x'\in V \setminus H$ in the same stratum as $x$), while on $\pi_1(S^1)$, it can be seen from (\ref{relinf}) that $lk_x$ acts by sending the generator (which coincides with the homotopy class of the meridian loop $\gamma_{\infty}$ about $H$) to ${-d}$.  The acyclicity at $x \in V \cap H$ then follows by the K\"unneth formula, since the homotopy factors of $\cU_x$, endowed with the corresponding homomorphisms and representations induced from the pair $(lk_x,\rho_x)$, are acyclic.
\end{proof}

\subsection{Sheaf (co)homology interpretation of twisted Alexander modules}\label{sh}
For the remaining of the paper, we will employ the language of perverse sheaves for relating local and global properties of twisted Alexander invariants. For this purpose, we first rephrase the definition of twisted Alexander modules as the (co)homology of a certain local system defined on the complement $\cU$.

Let $\cL$ be the local system of $\Ft$-modules on $\cU$, with stalk $\Ft \otimes_{\bF} \bV$, and action of the fundamental group $$\pi_1(\cU) \lra Aut(\Ft \otimes_{\bF} \bV)\cong GL_{\ell}(\Ft)$$ given by $$[\alpha] \mapsto t^{\eps(\alpha)} \otimes \rho(\alpha).$$
(Here $\ell$ denotes as before the rank of the representation $\rho$.) 
Then it is clear from the definition of the homological twisted Alexander modules that we have the following isomorphism of $\Ft$-modules:
\be\label{hsh}
\TA_i(\cU,\Ft) \cong H_i(\cU,\cL).
\ee
Note that $\Ft$ has a natural involution, denoted by $ \bar{\ \cdot \ } $ and defined by $t \mapsto t^{-1}$. Let $\cL^{\vee}$ be the local system on $\cU$ dual to $\cL$. Then $\cL^{\vee}\cong \bar{\cL}$, where $\bar{\cL}$ is the local system of $\Ft$-modules on $\cU$, with stalk $\Ft \otimes_{\bF} \bV^{\vee}$, but with the $\Ft$-module structure composed with the above involution. 
Then it follows as in \cite[pp.638]{KL} that the cohomological twisted Alexander modules of $(\cU,\eps,\rho)$ can be realized as:
\be\label{csh}
H^i_{\eps,\rho}(\cU,\Ft) \cong H^i(\cU,\cL^{\vee}).
\ee

If $x \in V$, let $i_x:\cU_x:=\cU \cap B_x \hookrightarrow \cU$ denote the inclusion of local complement at $x$, with corresponding induced local pair $(\eps_x,\rho_x)$ as in Section \ref{loca}.
Let $\cL_x:=i_x^*\cL$ be the restriction of the local system $\cL$ to $\cU_x$, i.e., $\cL_x$ is defined via the action of $(\eps_x,\rho_x)$. Then, for any $k \in \Z$,
 the local $k$-th (co)homological twisted Alexander modules at $x$ can be described as: $$H_k^{\eps_x,\rho_x}(\cU_x,\Ft)\cong H_k(\cU_x,\cL_x) \  \text{ and } \  H^k_{\eps_x,\rho_x}(\cU_x,\Ft)\cong H^k(\cU_x,\cL_x^{\vee}).$$


\subsection{Local-to-global analysis. Divisibility results}\label{lg}
In this section, we assume that the projective hypersurface $V$ is in general position at infinity. By Theorem \ref{t3} and Corollary \ref{cst}, the (co)homological twisted Alexander modules $H^{\eps,\rho}_{i}(\cU,\Ft)$, resp. $H_{\eps,\rho}^{i}(\cU,\Ft)$, are torsion $\Ft$-modules for any $0 \leq i \leq n$. Following Definition \ref{tap}, we denote by $\Delta^{\eps,\rho}_{i,\cU}(t)$, resp., $\Delta_{\eps,\rho,\cU}^{i}(t)$, with $0\leq i \leq n$, the corresponding twisted Alexander polynomials.

The sheaf theoretic realization of twisted Alexander modules in Section \ref{sh} allows us the use of perverse sheaves (or intersection homology) which, when coupled with homological algebra techniques, provide a concise relationship between the global twisted Alexander invariants of complex hypersurface complements and the corresponding local ones at singular points (respectively, at infinity). For simplicity of exposition, we choose to formulate our results in this section in cohomological terms, but see also Remark \ref{finrem} below. Our approach  is similar to  \cite[Section 3]{DM}.

\medskip

We work with sheaves of $\Ft$-modules. For a topological space $Y$, we denote by $D^b_c(Y;\Ft)$ the derived category of complexes of sheaves of $\Ft$-modules on $Y$ with constructible cohomology, and we let $\Perv(Y)$ be the abelian category of perverse sheaves of $\Ft$-modules on $Y$.

\medskip


The first result of this section singles out the contribution of the loop ``at infinity'' $\gamma_{\infty}$  to the global twisted Alexander invariants, and it can be regarded as a high-dimensional generalization (and for arbitrary singularities) of Corollary \ref{c1}, where $\gamma_{\infty}$ plays the role of $x_0^{-1}$ in loc.cit.:

\bt\label{t4} Let $V \subset \CP$ be a projective hypersurface in general position (with respect to the hyperplane $H$) at infinity, with complement $\cU=\CP\setminus (V \cup H)$. Fix a non-trivial epimorphism $\eps:\pi_1(\cU) \to \Z$ and a rank $\ell$ representation $\rho:\pi_1(\cU) \to GL(\bV)$. 
Then, for any $0 \leq i \leq n$, the zeros of the global cohomological Alexander polynomial $\Delta_{\eps,\rho,\cU}^{i}(t)$ are among those of the order of the cokernel of the endomorphism $t^{\eps(\gamma_{\infty})} \otimes \rho(\gamma_{\infty}) - Id \in \End(\Ft \otimes_{\bF} \bV)$.
\et

\begin{proof}
Let $\bC^{n+1}=\CP \setminus H$, and denote by
$u:\cU \hookrightarrow \bC^{n+1}$ and $v:
\bC^{n+1} \hookrightarrow \CP$ the two
inclusions. Since $\cU$ is smooth and $(n+1)$-dimensional, and $\cL^{\vee}$ is a local system on $\cU$, it follows that $\cL^{\vee}[n+1] \in \Perv(\cU)$.
Moreover, since $u$ is a quasi-finite affine
morphism, we also have that 
$$\cF^{\bullet}:=Ru_{\ast}(\mathcal{L}^{\vee}[n+1]) \in
\Perv(\bC^{n+1}),$$ e.g., see \cite[Theorem 6.0.4]{Sch}. But $\C^{n+1}$ is
an affine $(n+1)$-dimensional variety, so by Artin's vanishing
theorem for perverse sheaves (e.g., see \cite[Corollary 6.0.4]{Sch}), we obtain that:
\be\label{van1} \mathbb{H}^k(\bC^{n+1}, \cF^{\bullet})=0, \ \text{for all} \ k>0,\ee and
\be\label{van2} \mathbb{H}_c^k(\bC^{n+1}, \cF^{\bullet})=0, \ \text{for all} \ k<0.\ee
Let $a:\mathbb{CP}^{n+1} \to point$ be the constant map. Then:
\be\label{i1} \mathbb{H}^k(\bC^{n+1}, \cF^{\bullet})\cong H^{k+n+1}(\mathcal{U},
\mathcal{L}^{\vee})\cong H^k(Ra_{\ast}Rv_{\ast}\cF^{\bullet}).\ee Similarly, 
\be\label{i2} \mathbb{H}_c^k(\bC^{n+1}, \cF^{\bullet})\cong H^k(Ra_{!}Rv_{!}\cF^{\bullet}),\ee
where the last equality follows since $a$ is a proper map, hence $Ra_!=Ra_{\ast}$.

Consider the canonical morphism $Rv_{!}\cF^{\bullet} \to
Rv_{\ast}\cF^{\bullet}$, and extend it to the distinguished triangle:
\be\label{tri1} Rv_{!}\cF^{\bullet} \to Rv_{\ast}\cF^{\bullet} \to \mathcal{G}^{\bullet}
\overset{[1]}{\to}\ee in $D_c^b(\mathbb{CP}^{n+1};\Ft)$. Since
$v^{\ast}Rv_{!} \cong id \cong v^{\ast}Rv_{\ast}$, 
after applying $v^*$ to the above triangle we get that $v^*\cG\cong 0$, or equivalently, $\cG$ is supported on $H$.
Next, we apply $Ra_!=Ra_{\ast}$ to the distinguished triangle (\ref{tri1}) to obtain a new triangle in $D^b_c(point;\Ft)$:
\be\label{tri2}
Ra_{!}Rv_{!}\mathcal{F}^{\bullet}
\to Ra_{\ast}Rv_{\ast}\mathcal{F}^{\bullet} \to Ra_{\ast}\mathcal{G}^{\bullet}
\overset{[1]}{\to}
\ee
Upon applying the cohomology functor to the distinguished triangle (\ref{tri2}), and using the vanishing from (\ref{van1}) and (\ref{van2}) together with the identifications (\ref{i1}) and (\ref{i2}), we obtain that:
$$H^{k+n+1}(\mathcal{U},\mathcal{L}^{\vee})\cong \mathbb{H}^k(\mathbb{CP}^{n+1},
\mathcal{G}^{\bullet})\cong \mathbb{H}^k(H, \mathcal{G}^{\bullet}) \ \ \text{for} \
\ k < -1,$$ and $H^{n}(\mathcal{U}, \mathcal{L}^{\vee})$ is a
submodule of the $\Ft$-module $\mathbb{H}^{-1}(H, \mathcal{G}^{\bullet})$. 
So in order to prove the theorem, it remains to show that the $\Ft$-modules $\mathbb{H}^k(H, \mathcal{G}^{\bullet})$ are torsion for $k \leq -1$, and the zeros of their corresponding orders are amongs those of the order of the cokernel of $t^{\eps(\gamma_{\infty})} \otimes \rho(\gamma_{\infty}) - Id \in \End(\Ft \otimes_{\bF} \bV)\cong M_{\ell}(\Ft)$.

Note that  $\mathbb{H}^{k}(H,\mathcal{G}^{\bullet})$ is
the abutment of a hypercohomology spectral sequence with the $E_2$-term defined by
\be\label{ssg} E_2^{p,q}=H^p(H,\mathcal{H}^q(\mathcal{G}^{\bullet})).\ee This prompts us to investigate the stalk cohomology of $\cG^{\bullet}$ at points along $H$.

For $x \in H$, let us denote as before by $\cU_x=\cU \cap B_x$ the local complement at $x$, for $B_x$ a small ball in $\CP$ centered at $x$. Then we have the following identification:
\be\label{loc1} \mathcal{H}^q(\mathcal{G}^{\bullet})_x \cong
H^{q+n+1}(\mathcal{U}_x,\mathcal{L}^{\vee}_x),\ee
where $\cL_x$ is the restriction of $\cL$ to $\cU_x$. Indeed, the following isomorphisms of $\Ft$-modules hold:
\[
\begin{aligned}
 \mathcal{H}^q(\mathcal{G}^{\bullet})_x
& \cong \cH^q(Rv_*\cF^{\bullet})_x \\
& \cong \cH^{q+n+1}(Rv_*Ru_*\cL^{\vee})_x \\
& \cong \bH^{q+n+1}(B_x,R(v \circ u)_*\cL^{\vee}) \\
& \cong H^{q+n+1}(\mathcal{U}_x,\mathcal{L}^{\vee}_x).
\end{aligned}
\]

If $x \in H \setminus V$, then $\cU_x$ is homotopy equivalent to $S^1$, and the corresponding local system $\cL_x$ is defined by the action of $\gamma_{\infty}$, i.e., by multiplication by $t^{\eps(\gamma_{\infty})} \otimes \rho(\gamma_{\infty})$ on $\Ft \otimes_{\bF} \bV$. In particular, $H_*(U_x,\cL_x)$ is the homology of the complex of $\Ft$-modules:
$$0 \lra \Ft^{\ell} \xrightarrow{t^{\eps(\gamma_{\infty})} \otimes \rho(\gamma_{\infty}) - Id} \Ft^{\ell} \lra 0,$$ i.e., 
\be\label{loc2}
H_k(U_x,\cL_x)=
\begin{cases}
\cok(t^{\eps(\gamma_{\infty})} \otimes \rho(\gamma_{\infty}) - Id), & k=0 \\
0  , & k>0.
\end{cases}
\ee
If $x \in H \cap V$, then we know by Proposition \ref{ploc} that the local twisted Alexander modules $H^{k}(\mathcal{U}_x,\mathcal{L}^{\vee}_x)$ are $\Ft$-torsion, for all $k \in \Z$. Moreover, 
 in the notations of Proposition \ref{ploc}, $\cU_x \simeq \cU'_x \times S^1$, and the local system $\cL_x$ is an external tensor product, the second factor being defined by the action of $\gamma_{\infty}$ as in the previous case. So it follows from the K\"unneth formula that the zeros of the homological (hence also cohomological by the UCT) local twisted Alexander polynomials at points in $H \cap V$ are among those of the order of the cokernel of $t^{\eps(\gamma_{\infty})} \otimes \rho(\gamma_{\infty}) - Id \in \End(\Ft \otimes_{\bF} \bV)$.

By (\ref{loc1}) and the above calculations, it then follows that the $\Ft$-modules $\mathcal{H}^q(\mathcal{G}^{\bullet})_{x\in H}$ are torsion, and the zeros of their associated orders are among those of the order of the cokernel of $t^{\eps(\gamma_{\infty})} \otimes \rho(\gamma_{\infty}) - Id \in \End(\Ft \otimes_{\bF} \bV)$.
 Hence, by using the spectral sequence (\ref{ssg}), each hypercohomology group $\mathbb{H}^{k}(H,\mathcal{G}^{\bullet})$ ia a torsion $\Ft$-module, and the zeros of its associated order are among those of the order of the cokernel of $t^{\eps(\gamma_{\infty})} \otimes \rho(\gamma_{\infty}) - Id \in \End(\Ft \otimes_{\bF} \bV)$. This ends the proof of our theorem.
\end{proof}


\br If $\bF=\bC$ and $\eps=lk$ is the total linking number homomorphism, Theorem \ref{t4} implies that any root $\lambda$ of $\Delta^i_{\rho,\cU}(t)$, $i \leq n$, must satisfy the condition that $\lambda^d$ is an eigenvalue of $\rho(\gamma_{\infty})$, where $d=\sum_{i=1}^r n_i d_i$ is the degree of $V$.
If, in addition, 
 $\rho=triv$ is the trivial representation, the statement of Theorem \ref{t4} reduces to the fact that the zeros of the classical cohomological Alexander polynomials $\Delta^i_\cU(t)$, $i \leq n$, are roots of unity of order $d=\deg(V)$, a fact also shown in \cite{M06,DL,Liu} in the reduced case. 
\er

In the next theorem, we assume for simplicity of exposition that $V$ is a reduced hypersurface. 
Recall from Sections \ref{loca} and \ref{sh} that for any point $x$ in $V$, with local complement $\cU_x=\cU \cap B_x$, we get from $(\eps,\rho)$ an induced pair $(\eps_x,\rho_x)$ via the inclusion map $i_x:\cU_x \hookrightarrow \cU$. Moreover, the local twisted Alexander modules have a sheaf description in terms of the local system $\cL_x:=i_x^*\cL$ and its dual, namely, $H_k^{\eps_x,\rho_x}(\cU_x,\Ft)\cong H_k(\cU_x,\cL_x)$ and $H^k_{\eps_x,\rho_x}(\cU_x,\Ft)\cong H^k(\cU_x,\cL_x^{\vee})$, for all $k \in \Z$. We denote by $\Delta_{k,x}(t):=\Delta_{k,\cU_x}^{\eps_x, \rho_x}(t)$ and $\Delta^k_x(t):=\Delta^{k}_{\eps_x, \rho_x,\cU_x}(t)$ the {\it local} (co)homological twisted Alexander polynomials at $x$. 

Let us now assume also that $V$ is in general position at infity. Then if $x \in V \cap H$, in the notations of Proposition \ref{ploc} there is a homotopy equivalence $\cU_x \simeq \cU'_x \times S^1$, where $\cU'_x=B_x \setminus V$ and with the $S^1$-factor corresponding to the meridian loop about the hyperplane at infinity $H$. On the other hand, $\cU'_x$ is homeomorphic to any local complement $\cU_{x'}$ at a point $x'\in V \setminus H$ in the same stratum with $x$. So by the K\"unneth formula, the zeros of the local twisted Alexander polynomials $\Delta^k_x(t)$ of $(\cU_x,\eps_x,\rho_x)$ are among those associated to $(\cU_{x'},\eps_{x'},\rho_{x'})$, for $x' \in V \setminus H$ a nearby point in the same stratum of $V$ as $x$. For brevity, points of $V^a=V \setminus H$ will be referred to as {\it affine points of $V$}.

The next result shows that the zeros of the global twisted Alexander polynomials are completely determined by those of the local twisted Alexander polynomials at (affine) points along some irreducible component of $V$:

\bt\label{t5} Let $V \subset \CP$ be a reduced hypersurface in general position at infinity, with complement $\cU=\CP\setminus (V \cup H)$, and let $V_1$ be a fixed irreducible component of $V$. Fix a non-trivial epimorphism $\eps:\pi_1(\cU) \to \Z$, a rank $\ell$ representation $\rho:\pi_1(\cU) \to GL(\bV)$, and a non-negative integer $\sigma$. 
If $\lambda \in \bF$ is not a root of the $i$-th local twisted Alexander polynomial $\Delta^i_x(t)$ for any $i<n+1-\sigma$ and any (affine) point $x \in V_1 \setminus H$, then $\lambda$ is not a root of the global twisted Alexander polynomial $\Delta^i_{\eps,\rho,\cU}(t)$ for any $i<n+1-\sigma$.
\et

\begin{proof} 
First note that by the transversality assumption and K\"unneth, it follows by the above considerations that the hypothesis on local Alexander polynomials implies that $\lambda$ is not a root of the $i$-th local twisted Alexander polynomial $\Delta^i_x(t)$ for any $i<n+1-\sigma$ and any point $x \in V_1$ (including points in $V_1 \cap H$).

As in the proof of Theorem \ref{t4}, after replacing $\C^{n+1}$ by $\mathcal{U}_1=\mathbb{CP}^{n+1} \setminus V_1$, it follows that for $k \leq -1$, $H^{k+n+1}(\cU,\cL^{\vee})$ is a submodule of $\bH^k(\CP,\cG^{\bullet})$, where $\cG^{\bullet}$ is now a complex of sheaves of $\Ft$-modules supported on $V_1$. It thus suffices to show that $\bH^k(\CP,\cG^{\bullet})$, $k < -\sigma$, is a torsion $\Ft$-module whose order does not vanish at $\lambda$.

As in (\ref{loc1}), the cohomology stalks of $\cG^{\bullet}$ at any $x \in V_1$ are given by
 $$\mathcal{H}^q(\mathcal{G}^{\bullet})_x \cong H^{q+n+1}(\mathcal{U}_x,\mathcal{L}^{\vee}_x),$$ and these are all torsion $\Ft$-modules by Proposition \ref{ploc}.  Therefore, for a fixed $x \in V_1$ the fact that $\lambda$ is not a root of $\Delta^i_x(t)$ for any $i<n+1-\sigma$ is equivalent to the assertion that the order of  $\mathcal{H}^q(\mathcal{G}^{\bullet})_x$ does not vanish at $\lambda$ for all $i < -\sigma$. The desired claim follows now by using the hypercohomology spectral sequence with $E_2$-term defined by
$E_2^{p,q}=H^p(V_1,\mathcal{H}^q(\mathcal{G}^{\bullet}))$, which computes 
the groups $\mathbb{H}^{k}(V_1,\mathcal{G}^{\bullet})\cong \bH^k(\CP,\cG^{\bullet})$.
\end{proof}

\br Note that the proofs of Theorems \ref{t4} and \ref{t5} indicate that we can give a more general condition than transversality with respect to $H$ in order to conclude that the global cohomological twisted Alexander modules $H^i_{\eps,\rho}(\cU;\Ft)$ are torsion for all $i \leq n$. Indeed, it suffices to assume that the pair $(\eps,\rho)$ is acyclic along $V\cap H$ (or even $V_1 \cap H$, in the context of Theorem \ref{t5}). Of course this assumption is satisfied if $V$ is in general position at infinity, as Proposition \ref{ploc} shows. But there are other instances when it is satisfied, like in the examples discussed in Section \ref{ex}.
\er

\br\label{finrem}
Let us conclude with a few observations about other possible approaches for studying twisted Alexander-type invariants of hypersurface complements. 

 If $\bF=\bC$, one can argue as in \cite{DL} if similar divisibility results are desired for the homological twisted Alexander polynomials. In more detail, the study of such twisted homological invariants is reduced via a twisted version of the Milnor sequence to studying the vanishing (except in the middle degree) of the homology groups $H_k(\cU,\cL_{\lambda} \otimes \bV_{\rho})$ (or equivalently, of cohomology groups $H^k(\cU,\cL_{\lambda} \otimes \bV_{\rho})$), where $\cL_{\lambda}$ is the rank-one $\bC$-local system on $\cU$ defined by the character $$\pi_1(\cU) \xrightarrow{\varepsilon} \Z \xrightarrow{1\mapsto \lambda} \bC^*.$$ The language of $\C$-perverse sheaves can then be employed as in the proofs of Theorems \ref{t4} and \ref{t5} to get the desired vanishing, thus providing a twisted generalization of results from \cite{M06,DL}. 

Alternatively, one can use the approach from \cite{M06,LM} to study the (co)homological twisted Alexander invariants by using the associated {\it residue complex} $\cR^{\bullet}$ of $\cU$, which is defined as the cone of the natural morphism $Rj_!\cL \lra Rj_*\cL$, for $j:\cU \hookrightarrow \CP$ the inclusion map. 

Lastly, such results can also be derived by using more elementary techniques as follows: first,  by transversality and a Lefschetz-type argument one can reduce, as in \cite{L94}, the study of the twisted Alexander modules of $\cU$ to those of a regular neighborhood $\cN$ in $\bC^{n+1}$ of the affine part $V^a$ of $V$; secondly, Alexander-type invariants of $\cN$ can be computed via the Mayer-Vietoris spectral sequence for the induced stratification of such a neighborhood. 

We leave the details and precise formulations as an exercise for the interested reader.
\er







\end{document}